\documentclass[fontsize=12pt,a4paper,headings=normal,
twoside=false,leqno,parskip=half-,abstract=true]{scrartcl}
\usepackage[english]{babel}
\usepackage[utf8]{inputenc}
\usepackage{hyperref}
\hypersetup{
 pdftitle={Lyapunov function},
 pdfauthor={Phillipo Lappicy, Bernold Fiedler},
 colorlinks=true,
 linkcolor=blue,
 citecolor=blue,
 filecolor=blue,
 urlcolor=blue}

\usepackage{graphicx}
\usepackage[format=plain,labelfont=bf,font=small]{caption}
\usepackage{xcolor}
\usepackage[arrow, matrix, curve]{xy}
\usepackage{float}

\usepackage{caption}
\usepackage{tabulary}
\usepackage{array}
\newcolumntype{N}[1]{>{\centering\arraybackslash}m{#1}}

\usepackage{amsmath,amsthm}
\usepackage{amssymb} 

\makeatletter
\newcommand{\tpitchfork}{%
  \vbox{
    \baselineskip\z@skip
    \lineskip-.52ex
    \lineskiplimit\maxdimen
    \m@th
    \ialign{##\crcr\hidewidth\smash{$-$}\hidewidth\crcr$\pitchfork$\crcr}
  }%
}
\makeatother
\usepackage{latexsym}

\usepackage[notref,notcite,color,final 
]{showkeys}

\definecolor{refkey}{rgb}{1,0,0}
\definecolor{labelkey}{rgb}{1,0,0}

\usepackage{tikz}

\usepackage[textwidth=2cm,textsize=small,backgroundcolor=none]{todonotes}

  \mathchardef\ordinarycolon\mathcode`\:
  \mathcode`\:=\string"8000
  \begingroup \catcode`\:=\active
    \gdef:{\mathrel{\mathop\ordinarycolon}}
  \endgroup

\theoremstyle{plain}
\newtheorem{thm}{Theorem}[section]

\newtheorem{cor}[thm]{Corollary}

\begin{document}

\title{{\LARGE{A Lyapunov function for  fully nonlinear parabolic equations in one spatial variable}}}

\author{
 \\
{~}\\
Phillipo Lappicy* and Bernold Fiedler**\\
\vspace{2cm}}

\date{version of \today}
\maketitle
\thispagestyle{empty}

\vfill

$\ast$\\
Instituto de Ciências Matemáticas e de Computação\\
Universidade de S\~ao Paulo\\
Avenida trabalhador são-carlense 400\\
13566-590, São Carlos, SP, Brazil\\
\\
$\ast \ast$\\
Institut für Mathematik\\
Freie Universität Berlin\\
Arnimallee 3\\ 
14195 Berlin, Germany\\


\newpage
\pagestyle{plain}
\pagenumbering{arabic}
\setcounter{page}{1}

\begin{abstract}
Lyapunov functions are used in order to prove stability of equilibria, or to indicate a gradient-like structure of a dynamical system. Zelenyak (1968) and Matano (1988) constructed a Lyapunov function for quasilinear parabolic equations. We modify Matano's method to construct a Lyapunov function for fully nonlinear parabolic equations under Dirichlet and mixed nonlinear boundary conditions of Robin type.

\ 

\end{abstract}

\section{Main results}

\numberwithin{equation}{section}
\numberwithin{figure}{section}
\numberwithin{table}{section}

We consider the  scalar fully nonlinear, strictly parabolic partial differential equation
\begin{equation}\label{FULLYEQ}
    f(x,u,u_x,u_{xx},u_t)=0
\end{equation}
for $x\in (0,1)$. Here indices abbreviate partial derivatives. We assume $f$ satisfies the parabolicity condition 
\begin{equation}\label{par}     
    f_q\cdot f_r<0    
\end{equation}
for every argument $(x,u,p,q,r)=(x,u,u_x,u_{xx},u_t)$. 
For simplicity we assume $f \in C^1$ and, only at $r=0$, also $f_p,f_q,f_r \in C^1$.

We consider \eqref{FULLYEQ} under two types of \emph{separated boundary conditions} at $x=\iota\in \{ 0,1\}$. For each boundary $x=\iota$, separately, we either assume homogeneous Dirichlet boundary conditions,
\begin{equation}\label{Dir}     
    u=0,
\end{equation}
or nonlinear boundary conditions of Robin type
\begin{equation}\label{Rob}     
    u_x=b^\iota(u).
\end{equation}

We assume $b^\iota\in C^1$. Neumann boundary conditions appear as the special case $b^\iota(u)=0$. 
We observe \eqref{Rob} is equivalent to the seemingly more general case $b^\iota(u,u_x)=0$, provided that latter equation can be solved for $u_x$, uniquely and globally. 
See \cite{Amann88} and \cite{Lunardi95} for abstract settings involving nonlinear boundary conditions of the type \eqref{Rob}.

We consider \emph{classical solutions} $u=u(t,x)$ of \eqref{FULLYEQ} with separated boundary conditions \eqref{Dir} or \eqref{Rob}. This means that we assume boundedness and continuity of $u$ and its partial derivatives $u_t,u_{xx}$, for $t\in (0,T)$, $x\in (0,1)$. 
Initial conditions $u(t,x)=u_0(x)$, at $t=0$, are therefore approached continuously, for $t\searrow 0$, in the sup-norm with respect to $x\in (0,1)$. Similarly, the boundary values $u(t,x)$, $u_x(t,x)$ exist as continuous boundary limits $x\to \iota\in \{0,1\}$, and are assumed to either satisfy \eqref{Dir} or \eqref{Rob} at that boundary. 

Below we construct a \emph{Lyapunov energy function}
\begin{equation}\label{IntroLyap}     
    E:=\int_{0}^{1} L(x,u,u_x)\,dx\,, \qquad \text{ such that } \frac{dE}{dt}<0
\end{equation}
along nonequilibrium classical solutions $u=u(t,x)$ of \eqref{FULLYEQ}. 
Therefore the time dependent energy $t\mapsto E(u(t,.))$ decreases strictly, except at equilibria.

For quasilinear equations $f(x,u,p,q,r)=-r+a(x,u,p)q+h(x,u,p)$, a Lyapunov function $E$ was constructed, independently, by Zelenyak \cite{Zelenyak68} and Matano \cite{Matano88}. 
See also Hu \cite{Hu11}, and Fiedler et al \cite{FiedlerRagazzoRocha14}, for concise expositions of Matano's method. 
An analogous method for Jacobi systems, a spatially discrete variant, was developed in \cite{FiedlerGedeon99}. For an adaptation to diffusion with singular coefficients see \cite{LappicySing}. 

In order to construct a Lyapunov function $E$ as in \eqref{IntroLyap}, we rewrite the parabolic equation \eqref{FULLYEQ} suitably. 
Only two modifications of Matano's idea are required, along with proper assumptions.

First, we solve \eqref{FULLYEQ} for the diffusion variable $q=u_{xx}$ in terms of the other variables $(x,u,u_x,u_t)$. Indeed,  the parabolicity condition \eqref{par} implies $f_q\neq 0$.
This allows us to apply the implicit function theorem and to rewrite \eqref{FULLYEQ} as 
\begin{equation}\label{FULLYDIFF}
    u_{xx}=F(x,u,u_x,u_t).
\end{equation}
Here the parabolicity condition becomes
\begin{equation}\label{par2}     
    F_r>0    
\end{equation}
at any $(x,u,p,r)=(x,u,u_x,u_t)$, since implicit differentiation implies $F_r=-f_r/f_q>0$. We note that $F$ may not, and need not be, defined globally. 
We only consider $F(x,u,p,.)$ to be defined on a bounded or unbounded open interval of $r$, with limits $\pm \infty$ of $F$ at its boundaries.

Our second modification splits the function $F$ into two parts: one independent of $u_t$, and the other depending on $u_t$. 
In other words, we distinguish between the diffusion part $F^0$ related to the equilibrium ODE $u_t=0$, and the diffusion part $F^1$ related to time changing solutions. 
Specifically, to account for equilibria, at all, we define
\begin{equation}\label{F0}
    F^0(x,u,p):=F(x,u,p,0)   
\end{equation}
and assume $F^0$ is well-defined. 
(Otherwise we may put $F^0\equiv0$, artificially.
We do not comment further on this uninteresting case, below.)
To separate and highlight dependence on $u_t$,  we define
\begin{equation}\label{F1}
    F^1(x,u,p,r):=
    \begin{cases}
        (F(x,u,p,r)-F^0(x,u,p))/r\, &\text{ for } r\neq 0\\
        F_r(x,u,p,0)\, & \text{ for } r= 0.
    \end{cases}
\end{equation}
Our regularity assumption on $f$ then implies $F^1 \in C^1$ and $F_p^0 \in C^1$.

The parabolic equation \eqref{FULLYDIFF} can be rewritten as 
\begin{equation}\label{diffEQ}
    u_{xx}=F^0(x,u,u_x)+F^1(x,u,u_x,u_t)u_t\,.
\end{equation}
The parabolicity condition \eqref{par}, \eqref{par2} implies
\begin{equation}\label{par3}     
    F^1>0.
\end{equation}
Indeed, the monotonicity condition \eqref{par2} ensures \eqref{par3}, at $r=0$, as well as $\textrm{sign}(F-F^0) = \textrm{sign}(r)$, for $r\neq 0$. In the latter case, the numerator and denominator in \eqref{F1} have the same sign, yielding \eqref{par3}.

\begin{thm}\label{thm1}
    \emph{\textbf{Lyapunov Function}}
    
    Assume $C^1$-differentiability of $f, f_p, f_q, f_r$, as specified above, and parabolicity condition \eqref{par} all hold.
    
    Then there exists a Lagrange function $L=L(x,u,p)$, uniformly on bounded sets of $(u,p) \in\mathbb{R}^2 $, with $L, L_p, L_{pp}$ of class $C^1$, such that 
    \begin{equation}\label{Lyap}
        E:= \int_{0}^{1} L(x,u,u_x) dx
    \end{equation}
    is a Lyapunov function \eqref{IntroLyap} for the equation \eqref{FULLYEQ}. More precisely, bounded classical solutions $u(t,x)$ of \eqref{FULLYEQ} satisfy
    \begin{equation}
    \label{Lyapdecay}
        \frac{dE}{dt}=-\int_0^1 L_{pp}(x,u,u_x)F^1(x,u,u_x,u_t)\cdot  |u_t |^2dx
    \end{equation}
    with a strictly positive weight $L_{pp}F^1$.
\end{thm}

Semigroup settings of \eqref{FULLYEQ} on appropriate phase-spaces $X\supseteq C^1([0,1])$ have been provided by \cite{Lunardi95}, under additional assumptions. 
Even though fully nonlinear equations may not guarantee global existence of solutions, in general, some of them do. For instance, example 8.5.2 in \cite{Lunardi95} proves that $u_t=f(u_{xx})$, for $f\in C^3$ satisfying $f(0)=0$ and $f'(q)\geq \epsilon>0$, possesses global solutions.

Suppose that the resulting semiflow $u(t)$ is \emph{bounded and dissipative}, i.e. any solution $u(t)$ remains bounded for all times and enters some a priori fixed large ball in $X$, eventually. 
Suppose also that global orbits are precompact in $X$, that is, the closure of $\{ u(t) \text{ $ | $ } t\in\mathbb{R}\}$ is compact in $X$. Hence, there exists a \emph{global attractor} $\mathcal{A}\subset X$ as in Theorem 2.2 of \cite{Lady91}.
See also \cite{Hale88, ChVi02, CaLaRo13}.
As a consequence of the Lyapunov function \eqref{IntroLyap}, \eqref{Lyap}, the LaSalle invariance principle holds and bounded trajectories converge to (sets of) equilibria. 
See for example Henry \cite{Henry81}, Section 4.3. 
Next, let us suppose that equilibria of \eqref{FULLYEQ} are isolated, e.g.~due to hyperbolicity. 
By dissipativity, there are finitely many of them. 
Since the $\omega$-limit set is connected, it must consist of a single equilibrium. 
Similarly for the $\alpha$-limit set. Hence, the global attractor can be characterized as follows.

\begin{cor}\emph{\textbf{Attractor Decomposition}}

Suppose that the semiflow $u(t)$ of classical solutions of \eqref{FULLYEQ} in $X$ is bounded and dissipative with precompact orbits, such that all equilibria are hyperbolic. Then the global attractor $\mathcal{A}$ of \eqref{FULLYEQ} consists of finitely many equilibria, and of heteroclinics orbits $u(t)$ with pairs of distinct equilibria as $\alpha$- and $\omega$-limit sets.
\end{cor}

We conclude with the following open question: can one construct the global attractor $\mathcal{A}$ of \eqref{FULLYEQ} explicitly, based on a characterization by a permutation of boundary values or boundary slopes of equilibria, analogously to the semilinear cases \cite{FiedlerBrunovsky89}, \cite{FiedlerRocha96} and the quasilinear case \cite{Lappicy18}? This conjecture was already stated in Fiedler \cite{Fiedler96}, and the Lyapunov function presented here may serve as a first step towards proving it.

Towards that goal, however, one still needs to establish zero numbers, transversality of stable and unstable manifolds, construct shooting curves, study their Sturm permutations, calculate Morse indices, and prove liberalism. We are optimistic that these difficulties will be overcome, in due time, along the promising lines of \cite{Lappicy18}.

\textbf{Acknowledgment.} We are indebted to Marek Fila for helpful discussions concerning nonlinear boundary conditions. Phillipo Lappicy was supported by FAPESP, Brasil, grant number 2017/07882-0. Bernold Fiedler was partially supported by SFB 910 of the Deutsche Forschungsgemeinschaft, and some generous libations of Cacha\c{c}a de Jamb\'u.

\section{Proof}

As a prerequisite, we first comment on higher regularity of classical solutions $u=u(t,x)$ of the nonlinear PDE \eqref{FULLYEQ} with Dirichlet or nonlinear Robin boundary conditions \eqref{Dir} or \eqref{Rob}, respectively. We recall our standing assumption of continuity and boundedness of $u$ and the partial derivatives $q=u_{xx}$, $r=u_t$ for $t\in (0,T)$, $x\in (0,1)$. 

Differentiation of \eqref{FULLYEQ} with respect to $x$, along any classical solution $u=u(t,x)$, provides an inhomogeneous linear nonautonomous parabolic PDE for $p=u_x$, of the form
\begin{equation}\label{1x}
    p_t=a_0(t,x)p_{xx}+a_1(t,x)p_{x}+a_2(t,x)p+a_3(t,x).
\end{equation}
The bounded continuous coefficients $a_j(t,x)$ are given by partial derivatives of $f$, evaluated along the given solution $u(t,x)$.  
No matter whether the solution $u$ satisfies Dirichlet \eqref{Dir} or nonlinear Robin boundary conditions \eqref{Rob}, the resulting boundary conditions for $p=u_x$ are explicit, of inhomogeneous Dirichlet type, in terms of the given solution derivatives $u_x(t,x)$ at $x=\iota \in \{ 0,1\}$ and for $t>0$. 
Standard bootstrapping provides classical solutions $p$, starting at any positive initial time $t=t_0>0$. 

Analogously, differentiation of \eqref{FULLYEQ} with respect to $t$ provides a linear nonautonomous parabolic PDE for $r=u_t$, of the form
\begin{equation}\label{1t}
    r_t=a_0(t,x)r_{xx}+a_1(t,x)r_{x}+a_2(t,x)r
\end{equation}
with some other bounded continuous coefficients $a_j(t,x)$, derived as before. 
Only under Dirichlet boundary conditions \eqref{Dir} for $u$, however, we obtain Dirichlet boundary conditions $r=u_t=0$ for $r$. Nonlinear boundary conditions \eqref{Rob} for $u$, in contrast, linearize to the nonautonomous homogeneous form $r_x=b^\iota(t)r$, where $b^\iota(t):=b^\iota_u(u(t,\iota))$ remains bounded. Again, standard bootstrapping provides classical solutions $r$, from any positive initial time $t=t_0>0$ onwards.

After these preparations, we can embark on Matano's method with our two minor adaptations. 
We rewrite \eqref{FULLYEQ} as an equation \eqref{FULLYDIFF} for the diffusion $u_{xx}$ in terms of the other variables $(x,u,p,r)=(x,u,u_x,u_t)$. 
We then split the diffusion into two parts: one independent of $u_t$, and the other depending on $u_t$. See \eqref{F0}, \eqref{F1}.
This yields
\begin{equation}\label{DIFFEQ}
    u_{xx}=F^0(x,u,u_x)+F^1(x,u,u_x,u_t)u_t
\end{equation}
with parabolicity condition $F^1>0$; see \eqref{par3}.

We now enter Matano's proof. 
Let $p:=u_x$ and differentiate the definition \eqref{Lyap} of the Lyapunov function $E$ with respect to time $t$ along classical solutions $u(t,x)$ of \eqref{FULLYEQ},
\begin{equation}\label{part1}
    \frac{dE}{dt}= \int_{0}^1 \left( L_u u_t + L_p u_{x t}\right) dx .
\end{equation}
Here we used that $u_{xt}=p_t$ for the classical solution $p$ of \eqref{1x}. 
The Lagrange function $L$ depends on $(x,u,p)=(x,u,u_x)$, only.
It remains to determine $L$ such that $dE/dt<0$, except at equilibria. 
See \eqref{IntroLyap}.
Integrating the second term in \eqref{part1} by parts, we obtain
\begin{align}\label{part2}
\begin{split}
    \frac{dE}{dt}&= L_pu_t\Big|_0^1 +\int_{0}^1 \left( L_u -\frac{d}{dx}L_p  \right) u_t dx\\
    &=L_pu_t\Big|_0^1+ \int_{0}^1 \left( L_u -L_{px} -L_{pu}u_x -L_{pp}u_{xx}\right) u_t dx.
\end{split}
\end{align}
Indeed, continuity of $u,\, p=u_x$, and $r=u_t$, up to the boundary, holds for the classical solutions $u$ of \eqref{FULLYEQ}, $p$ of \eqref{1x}, and $r$ of \eqref{1t}. 
We carry out the differentiation of $L_p$ with respect to $x$ in \eqref{part2}, and substitute the splitting \eqref{F0} and \eqref{F1} of the nonlinearity $F$, as in \eqref{diffEQ}, to obtain
\begin{equation}\label{part2b}
        \frac{dE}{dt}=L_pu_t\Big|_0^1 +\int_{0}^1 \left( L_u -L_{px} -L_{pu}u_x -L_{pp}F^0\right) u_t dx-\int_{0}^1 L_{pp}F^1 \cdot |u_t|^2 dx.
\end{equation}
Again we used the classical, and hence bounded, solution $r=u_t$ of \eqref{1t}, to assure continuity and boundedness of the integrand, and hence existence of the integrals.

We seek to construct the Lagrange function $L$ such that the boundary terms vanish, the parenthesis in the first integral \eqref{part2b} also vanishes, and $L_{pp}>0$. This yields a Lyapunov function such that
\begin{equation}\label{Lyaput2}
        \frac{dE}{dt}=-\int_{0}^1 L_{pp}F^1 \cdot |u_t|^2 dx.
\end{equation}
Note $L_{pp}F^1> 0$, due to our requirement that $L_{pp}>0$, and the parabolicity condition $F^1>0$ of \eqref{par3}. 
To prove the theorem, it therefore remains to guarantee that there exists a function $L=L(x,u,p)$ satisfying $L_{pp}>0$ such that
\begin{equation}\label{Lyapquasi0}
    L_u-L_{px}-pL_{pu} -F^0L_{pp}=0
\end{equation}
for all $(x,u,p)\in (0,1)\times \mathbb{R}^2$,  and such that $L_pu_t=0$ on the boundary.
Although this part follows Matano \cite{Matano88}, almost verbatim, we include the necessary details. 
In this part, $(u,p)\in \mathbb{R}^2$ are real variables, quite simply, rather than solutions $u,u_x$ of PDEs, viz. functions of $(t,x)$.

Differentiating \eqref{Lyapquasi0} with respect to $p$, the terms $L_{pu}$ cancel:
\begin{equation}\label{Lyapquasi0_p}
    L_{ppx}+pL_{ppu}+F^0L_{ppp} = -F^0_pL_{pp} \,.
\end{equation}
To ensure $L_{pp}>0$, Matano makes an Ansatz 
\begin{equation}\label{Lpp}
    L_{pp}=\exp(g)>0 \,.
\end{equation}
Rewriting \eqref{Lyapquasi0_p} in terms of $g$ provides the linear first order PDE
\begin{equation}\label{Lyapquasig}
    g_{x}+pg_{u}+F^0g_{p} =-F^0_p \,.
\end{equation}
This can be solved by the method of characteristics: along the solutions of the auxiliary ordinary differential equation 
\begin{equation}
\label{quasishootlyap}
\begin{aligned}
    \dot{x}&=1\\
    \dot{u}&= p\\
    \dot{p}&=F^0(x,u,p)
\end{aligned}
\end{equation}
the function $g$ must satisfy
\begin{equation}\label{Lyapg}
    \dot{g}=-F^0_p(x,u,p),
\end{equation}
e.g. with the initial condition $g(0,u,p)=0$.

Note that the characteristic ODE \eqref{quasishootlyap} coincides with the equilibrium equation $u_{xx}=F^0(x,u,p)$, since $F^0=F$ for $u_t=0$. 
Without further assumptions on the nonlinearity $f$ in \eqref{FULLYEQ}, solutions to \eqref{quasishootlyap} may not exist on the whole required interval $x\in (0,1)$. 
The Lyapunov claim of Theorem \ref{thm1}, however, was only asserted to hold on any bounded subset of $(u,p)\in \mathbb{R}^2$. 
We may therefore use a cut-off for $F^0$, in \eqref{quasishootlyap}, to assert global existence of solutions to the characteristic equation.
Our differentiability assumptions on $f$ imply $g \in C^1$.

After this construction we now have to reverse gear and ascend from a $C^1$-function $g$ satisfying \eqref{Lyapquasig} to a Lagrange function $L$ satisfying \eqref{Lyapquasi0}. The general solution $L$ of $L_{pp}=\exp(g)$ can be obtained by integrating twice with respect to $p$:
\begin{equation}\label{L}
    L(x,u,p):=\int_0^p \int_0^{p_1} \exp(g(x,u,p_2))dp_2dp_1 + L^0(x,u)+L^1(x,u)p.
\end{equation}
This solves \eqref{Lyapquasi0_p}, by $C^1$-differentiability of $g$ with respect to to $(x,u,p)$. 
To ensure that $L$ is also a solution of \eqref{Lyapquasi0}, we have to determine $C^1$ integration ``constants'' $L^0$ and $L^1$, appropriately. 
Note $L, L_p, L_{pp} \in C^1$, as claimed.

Recall that \eqref{Lyapquasi0_p} was obtained through differentiation of \eqref{Lyapquasi0} with respect to $p$. Conversely, the left-hand side of \eqref{Lyapquasi0} is therefore independent of $p$. Hence \eqref{Lyapquasi0} is satisfied for all $p$ if it holds for $p=0$. At $p=0$, the construction of $L$ yields $L_p=L^1$, $L_{px}=L^1_x$ and $L_u=L^0_u$. Insertion in \eqref{Lyapquasi0} at $p=0$ yields 
\begin{equation}\label{FindingG0}
    L^0_u=L^1_x+F^0\exp(g)\,.
\end{equation}
Integrating with respect to $u$, we obtain
\begin{equation}
    L^0(x,u):=\int_0^{u} \left[L^1_x(x,u_1)+\mathrm{exp}(g(x,u_1,0))F^0(x,u_1,0)\right]du_1+L^{00}(x)\,.
\end{equation}
Of course we may omit $L^{00}(x)$, which in \eqref{Lyap} just integrates to an irrelevant additive constant for the Lyapunov function $E$.

To complete the proof it only remains to show that $L_pu_t$ vanishes at the boundaries $x=0,1$. 
At any boundary of Dirichlet type \eqref{Dir} this is trivial because $r=u_t=0$. 

In the case of a nonlinear Robin boundary condition \eqref{Rob} at only one of the boundaries, either at $x=\iota=0$ or at $1$, we have to choose $L$ such that $L_p(\iota,u,b^\iota(u))=0$. 
By our construction \eqref{L} of $L$, this is equivalent to 
\begin{equation}\label{G1}
    L^1(\iota,u):=-\int_0^{b^\iota(u)} \mathrm{exp}(g(\iota,u,p))dp \,,
\end{equation}
and we may choose $L^1$ to be independent of $x$.

If nonlinear Robin boundary conditions \eqref{Rob} are imposed at both boundaries, $x=\iota=0$ and $1$, then we define $L^1(\iota,u)$ as in \eqref{G1}. 
Linear interpolation $L^1(x,u):=(1-x)L^1(0,u)+xL^1(1,u)$ then provides  $L^1 \in C^1$ such that $L_p(\iota,u,b^\iota(u))=0$. In either case, this construction proves the theorem.

\section{Remarks}

We conclude with a few comments on modifications and generalizations of our result.

For nonlinearities of quasilinear type, i.e. 
\begin{equation}
    u_t=a(x,u,u_x)u_{xx}+h(x,u,u_x)
\end{equation}
with parabolicity condition $a(x,u,u_x)\geq \epsilon >0$, our definition of $F$ yields 
\begin{equation}
    u_{xx}=F(x,u,u_x,u_t):=(u_t-h)/a.
\end{equation} 
Hence $F^0:=-h/a$ and $F^1:=1/a$. In particular, our method recovers the Lyapunov function of Matano \cite{Matano88}, where
\begin{equation}
    \frac{dE}{dt} = -\int_{0}^1 \frac{L_{pp}}{a} \,  |u_t|^2 \, dx.
\end{equation}
Here $p:=u_x$, and $L$ satisfies the convexity condition $L_{pp}>0$. 
The standard semilinear case $a\equiv 1$ where $h$ depends only on $(x,u)$, is recovered via $g\equiv0$ and $L_{pp}\equiv1$.

Alternatively, the proof of the main theorem can also be obtained with the modified functions
\begin{equation}
\label{F0F1_alt}
\begin{aligned}
    F^0(x,u,p)&:=F(x,u,p,0)\,,\\   
    F^1(x,u,p,r)&:=F(x,u,p,r)-F^0(x,u,p)\,.
\end{aligned}
\end{equation}
Here $rF^1>0$ due to the parabolicity condition $F_r>0$. This yields a Lyapunov function such that
\begin{equation}
\label{Lyapdecay_alt}
    \frac{dE}{dt}:= -\int_{0}^1 L_{pp}F^1  u_t dx\,.
\end{equation}
Here $L_{pp}F^1u_t>0$, except at equilibria. We prefer the splitting \eqref{F0}, \eqref{F1} of the function $F$ over \eqref{F0F1_alt}, for purely aesthetical reasons, to extract the $L^2$ gradient flow decay term $|u_t|^2$ in \eqref{Lyapdecay}, explicitly, compared to \eqref{Lyapdecay_alt}.


For direct variational methods, standard Lyapunov, or energy, functions are often thought to be positive semidefinite, or at least bounded from below. 
For compact global attractors $\mathcal{A}$, in particular, the cut-off in the equation \eqref{quasishootlyap} of characteristics is justified: our Lagrange function $L$, and hence our Lyapunov function $E$, remains unchanged on $\mathcal{A}$ .
In more general settings, including solutions which blow up in finite time, boundedness of $E$ from below might of course fail. 
In fact, applications to fully nonlinear blow-up may require delicate discussions of the equilibrium characteristic equation \eqref{quasishootlyap}, beyond our crude cut-off.

\medskip


\begin{thebibliography}{10}

\small{
\bibitem{Amann88}
H.~Amann.
Parabolic evolution equations and nonlinear boundary conditions.
\emph{J. Diff. Eq.} \textbf{72}, 201 -- 269, (1988).

\bibitem{CaLaRo13}
A.N.~Carvalho, J.L.~Langa, J.C.~Robinson.
\emph{Attractors For Infinite-Dimensional Non-Autonomous Dynamical Systems.} 
Appl. Math. Sc \textbf{182}, Springer-Verlag New York, (2013).

\bibitem{ChVi02}
V.V.~Chepyzhov and M.I.~Vishik.
\emph{Attractors for {E}quations of {M}athematical {P}hysics}.
Colloq. AMS, Providence, (2002).

\bibitem{FiedlerBrunovsky89}
P.~Brunovsk{\'y} and B.~Fiedler.
Connecting orbits in scalar reaction diffusion equations II: The
  complete solution.
\emph{J. Diff. Eq.} \textbf{81}, 106--135, (1989).

\bibitem{Fiedler96}
B.~Fiedler.
Do global attractors depend on boundary conditions?
\emph{Doc. Math. J. DMV} \textbf{1}, 215--228, (1996).

\bibitem{FiedlerRagazzoRocha14}
B.~Fiedler, C.~Grotta-Ragazzo, and C.~Rocha.
An explicit Lyapunov function for reflection symmetric parabolic
  partial differential equations on the circle.
\emph{Russ. Math. Surv.} \textbf{69}, 27--42, (2014).

\bibitem{FiedlerGedeon99}
B.~Fiedler and T.~Gedeon.
A Lyapunov Function for Tridiagonal Competitive-cooperative Systems.
\emph{SIAM J. Math. Anal.} \textbf{30(3)}, 469--478, (1999).

\bibitem{FiedlerRocha96}
B.~Fiedler and C.~Rocha.
Heteroclinic orbits of semilinear parabolic equations.
\emph{J. Diff. Eq.} \textbf{125}, 239--281, (1996).

\bibitem{Hale88}
J.K.~Hale.
 \emph{Asymptotic Behavior of Dissipative Systems}.
 Math. Surv. \textbf{25}. AMS, Providence, (1988).


\bibitem{Henry81}
D.~Henry.
\emph {Geometric Theory of Semilinear Parabolic Equations}.
Springer-Verlag New York, (1981).

\bibitem{Hu11}
B.~Hu.
\emph{Blow-Up Theories for Semilinear Parabolic Equations}.
Springer-Verlag Berlin, (2011).

\bibitem{Lady91}
O.~Ladyzhenskaya.
\emph{Attractors for Semi-groups and Evolution Equations}.
Cambridge University Press, (1991).

\bibitem{Lappicy18}
P.~Lappicy.
Sturm attractors for quasilinear parabolic equations.
\emph{J. Diff. Eq.}, \textbf{265}, 4642--4660, (2018).

\bibitem{LappicySing}
P.~Lappicy.
Sturm attractors for quasilinear parabolic equations with singular coefficients.
\emph{arXiv:1806.04019}, (2018).

\bibitem{Lunardi95}
A.~Lunardi.
\emph{Analytic Semigroups and Optimal Regularity in Parabolic Problems}.
Springer Basel, (1995).

\bibitem{Matano88}
H.~Matano.
Asymptotic behavior of solutions of semilinear heat equations on
  $S^1$.
\emph{Nonlinear Diffusion Equations and Their Equilibrium States II,}
  eds. W.-M. Ni, L. A. Peletier, J. Serrin, 139--162, (1988).

\bibitem{Zelenyak68}
T.~I. Zelenyak.
Stabilization of solutions of boundary value problems for a second
  order parabolic equation with one space variable.
 \emph{Differ. Uravn.} \textbf{4}, 34--45, (1968).
 
 }
\end{thebibliography}

\end{document}